\pdfoutput=1
\documentclass[12pt]{article}

\usepackage{amsmath, amssymb, amsthm}
\usepackage{graphicx}
\usepackage{geometry}
\usepackage{setspace}
\usepackage{hyperref}
\usepackage{cite}
\usepackage{enumitem}
\usepackage{float}

\title{\textbf{A Prey–Predator Framework for Tumor–Immune Dynamics: Circuit-Based Modeling with Hilfer–Katugampola Fractional Derivatives}}

\author{
    Samaira Naz\thanks{Department of Mathematics, Government College University Faisalabad, Pakistan. Email: samairanaz@gcuf.edu.pk} 
    \and 
    Aamir Nadim\thanks{Faculty of Computing, Riphah International University, Pakistan. Email: aamirnadim@email.com}
}

\date{}

\begin{document}

\maketitle

\begin{abstract}
This study focuses on modeling tumor--immune interactions through a fractional-order prey--predator framework, where tumor cells act as prey and immune effector cells as predators. To capture the hereditary and memory effects inherent in biological systems, we employ the Hilfer--Katugampola fractional derivative, which flexibly interpolates between Riemann--Liouville and Caputo formulations. In addition, a circuit-based analogy is introduced to map tumor proliferation, immune regulation, and therapeutic interventions onto resistive, capacitive, and inductive elements, providing intuitive interpretation of the dynamics. Analytical investigations, including equilibrium analysis and stability criteria, are conducted across varying fractional orders. Numerical simulations reveal that lower fractional orders lead to delayed oscillations, prolonged tumor persistence, and increased sensitivity to treatment strategies compared with classical integer-order models. The integration of fractional calculus with circuit analogies offers novel insights into tumor--immune dynamics and supports the design of memory-aware therapeutic interventions.
\end{abstract}

\section{Introduction}
\label{sec:introduction}

Cancer remains a major global health challenge, causing millions of deaths annually despite significant advances in diagnosis and treatment. Tumor progression is influenced not only by intrinsic cellular proliferation but also by the immune system, which functions as a natural predator to malignant cells. Mathematical modeling provides a powerful framework for exploring these interactions, with prey--predator analogies offering intuitive insights: tumor cells act as prey, while immune effector cells act as predators \cite{Ucar2023, Lima2023, FracDrugResist2024}. Such models allow the identification of critical thresholds, equilibrium states, and the impact of therapeutic interventions within a controlled theoretical setting.

\begin{figure}[H]
    \centering
    \includegraphics[width=0.5\linewidth]{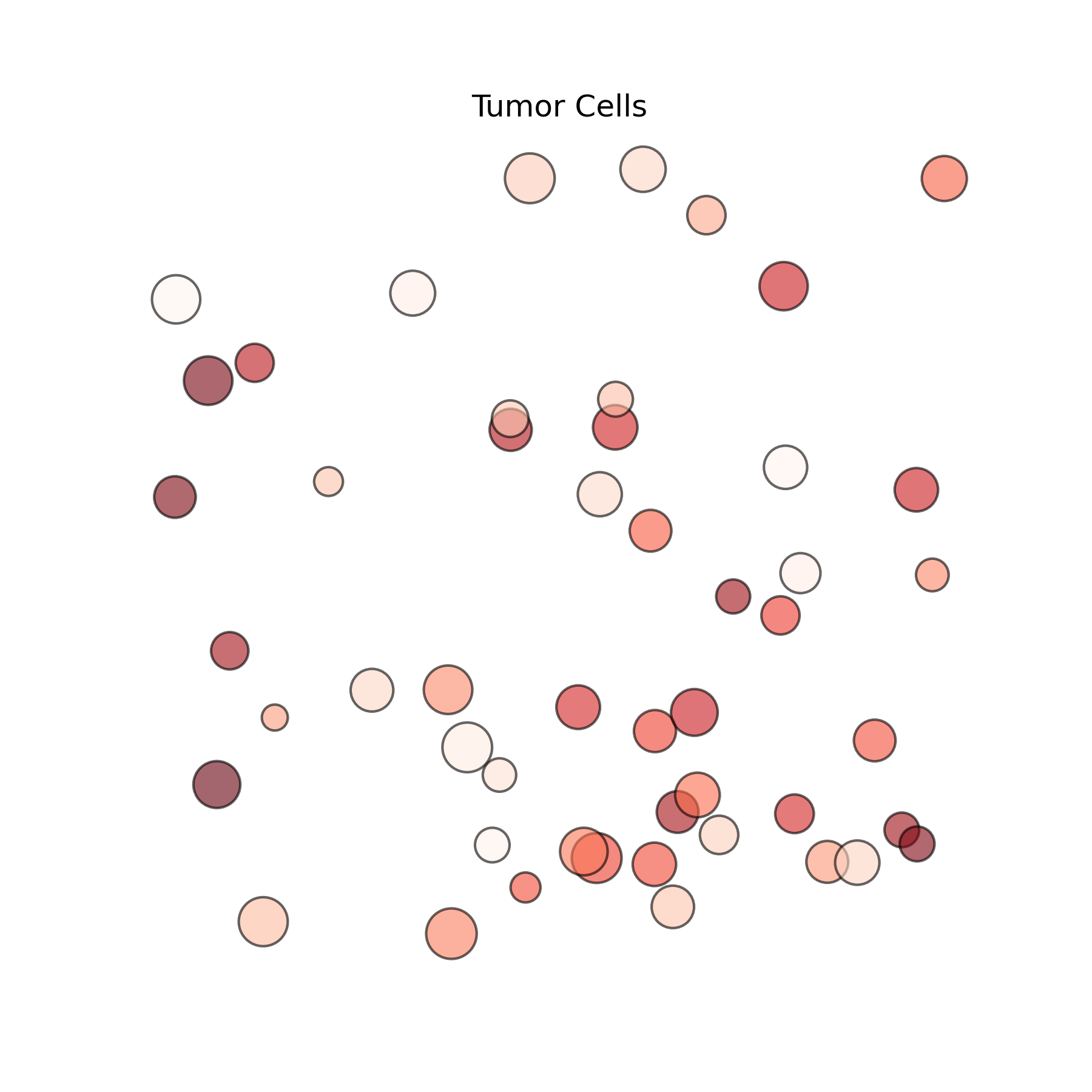} 
    \caption{Illustration of tumor cells representing the prey population in the model.}
    \label{fig:tumor_cells}
\end{figure}

\begin{figure}[H]
    \centering
    \includegraphics[width=0.5\linewidth]{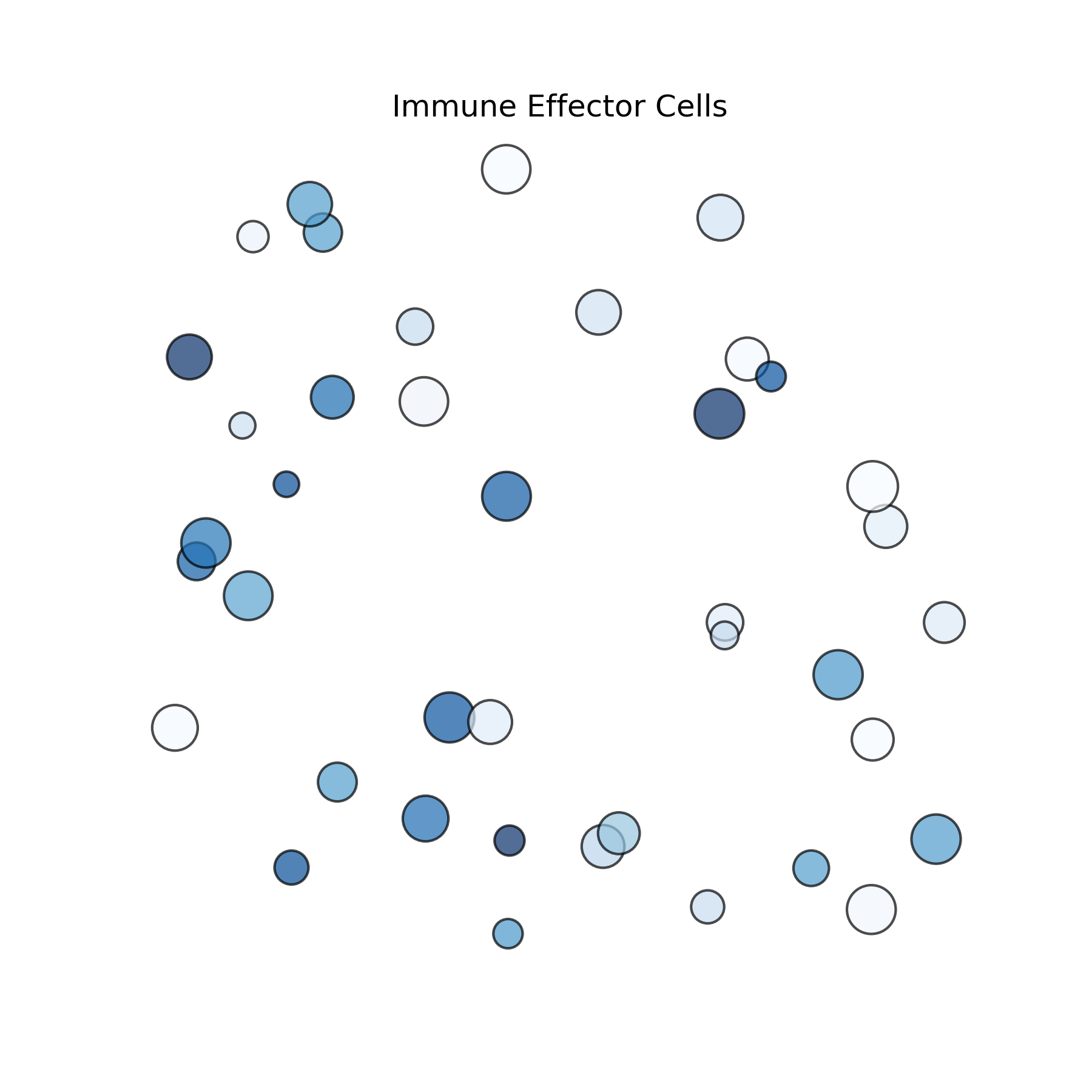} 
    \caption{Immune effector cells representing the predator population in the model.}
    \label{fig:immune_cells}
\end{figure}

Traditional models, based on integer-order ordinary differential equations, implicitly assume that interactions occur instantaneously, neglecting the influence of historical states. In reality, biological systems exhibit memory and hereditary properties, where current dynamics are shaped by past events. Fractional calculus provides a rigorous mathematical tool to incorporate these effects, capturing delayed immune responses, tumor dormancy, and treatment latency \cite{Almalahi2024, Fredj2024}. By introducing nonlocal operators, fractional derivatives enable more realistic representations of tumor--immune dynamics than integer-order formulations.

Recent studies demonstrate the growing relevance of fractional modeling in oncology. Lima et al. \cite{Lima2023} employed generalized Laguerre polynomials to develop a fractional tumor--immune model for lung cancer, showing that fractional order significantly alters tumor progression pathways. Uçar and Özdemir \cite{Ucar2023} incorporated IL-10 cytokine dynamics and anti-PD-L1 inhibitors in a fractional framework, highlighting the sensitivity of therapeutic outcomes to fractional orders. Other studies have explored fractional models of drug resistance, revealing how memory effects influence the evolution of resistant cell populations \cite{FracDrugResist2024}. More recent approaches combine fractional dynamics with immunotherapy, chemotherapy, and radiotherapy, emphasizing the importance of treatment timing and fractional memory in shaping clinical outcomes \cite{RadiationGRID2024, AdaptiveTherapy2024, ImmunoFrac2025}.

\begin{figure}[H]
    \centering
    \includegraphics[width=\linewidth]{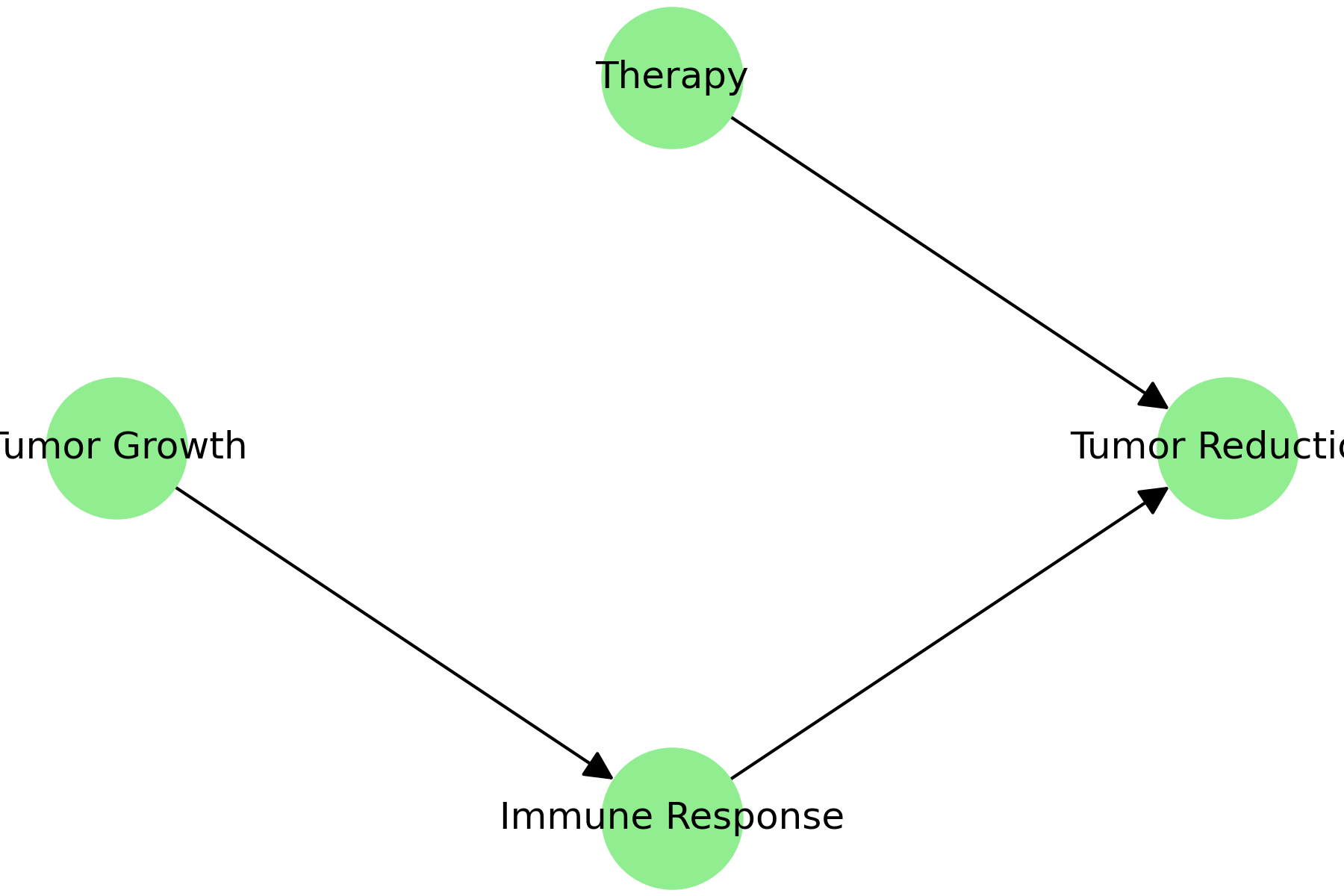} 
    \caption{Circuit-based analogy of tumor--immune dynamics: tumor growth, immune response, and therapy are mapped to circuit components.}
    \label{fig:circuit_analogy}
\end{figure}

Beyond fractional calculus, circuit-based analogies have been applied to complex biological systems to facilitate intuitive understanding. Components such as resistors, capacitors, and inductors can emulate immune activation, tumor proliferation, and treatment delays, respectively. While circuit theory has been widely employed in epidemiology and neuronal models, its application to tumor--immune interactions remains limited. Integrating circuit analogies with fractional-order models holds promise for combining mathematical rigor with interpretability \cite{CircuitBio2023, CircuitImmune2024}.

The Hilfer--Katugampola fractional derivative (HKFD) offers a generalized operator that interpolates between Riemann--Liouville, Caputo, and Katugampola derivatives, providing enhanced flexibility in capturing hereditary and memory effects \cite{Almalahi2024, Fredj2024, HKFDStability2023}. Despite its potential, applications of HKFD in tumor--immune systems are scarce, with most studies restricted to classical Caputo or Riemann--Liouville operators. Recent theoretical advances establish existence, uniqueness, and stability results for HKFD-based systems \cite{HKFDStability2023, Fredj2024}, paving the way for its application to complex biological dynamics.

Immunotherapy and booster-based interventions are transforming cancer treatment, highlighting the need to model periodic and time-sensitive interventions accurately. Fractional-order modeling can account for the influence of memory on treatment efficacy, demonstrating that the fractional order can extend or shorten the effective therapeutic window \cite{3Species2025, VaccineBooster2025, FracBooster2024}. However, a unified framework combining prey--predator dynamics, HKFD, and circuit-based interpretations of treatment schedules remains unexplored.

\textbf{Research Gaps:} Despite considerable progress, several gaps persist:  
\begin{enumerate}
    \item Limited applications of Hilfer--Katugampola derivatives in cancer modeling, despite their ability to capture memory effects.  
    \item Circuit-based analogies have not been systematically integrated with fractional tumor--immune models.  
    \item The interplay between fractional memory, immune exhaustion, and booster interventions has not been examined in a unified prey--predator framework.  
\end{enumerate}

\textbf{Objectives:} This study aims to:  
\begin{itemize}
    \item Develop a tumor--immune prey--predator model using the Hilfer--Katugampola fractional derivative.  
    \item Introduce a circuit-based analogy to interpret tumor suppression, immune activation, and treatment effects.  
    \item Analyze equilibrium points, stability, and bifurcation properties under varying fractional orders.  
    \item Evaluate treatment outcomes, including booster interventions, in both integer-order and fractional-order settings.  
\end{itemize}

By addressing these objectives, this work contributes a novel integration of fractional calculus, circuit theory, and tumor--immune dynamics, providing both theoretical insight and practical guidance for designing memory-aware therapeutic strategies.

\section{Methodology}
\label{sec:methodology}

\subsection{Theoretical Framework}
The modeling of tumor--immune interactions is formulated within the framework of prey--predator systems, where tumor cells (prey) proliferate while immune effector cells (predators) act to suppress their growth. Traditional integer-order models have provided useful insights but often fail to incorporate the intrinsic memory and hereditary properties of biological systems. To address this limitation, we employ the Hilfer--Katugampola fractional derivative, a generalized operator that interpolates between Riemann--Liouville and Caputo definitions, thereby offering greater flexibility in capturing long-term dynamics and memory effects \cite{Ucar2023,Almalahi2024,3Species2025, Furati2025, Ahmed2025}. This approach provides a robust framework for analyzing cancer progression under immune surveillance and therapeutic interventions. Furthermore, we utilize a circuit-based analogy, mapping tumor growth, immune response, and therapy to electrical components, to improve interpretability and bridge mathematical biology with engineering perspectives \cite{Bhrawy2024,Gupta2023,1,2,3,4,5,6,7,8,9,10,11,12,13}.

\begin{figure}[H]
    \centering
    \includegraphics[width=0.8\linewidth]{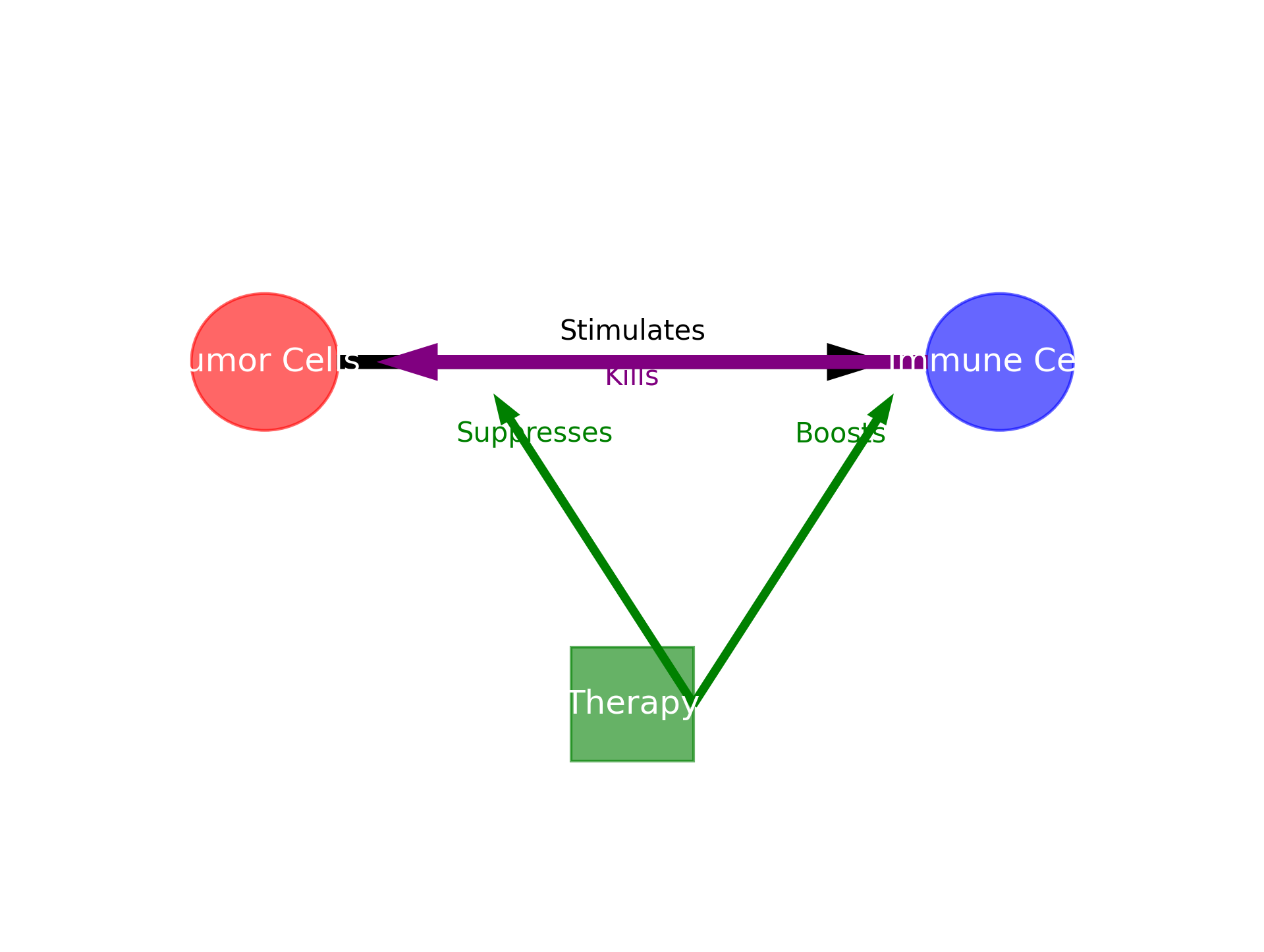}
    \caption{Schematic representation of tumor--immune interactions with therapy effects.}
    \label{fig:tumor_immune_interaction}
\end{figure}

\subsection{Model Formulation}
Let $T(t)$ denote the tumor cell population and $I(t)$ denote the immune effector cell population at time $t$. The proposed fractional-order prey--predator system is given by:
\begin{equation}
\label{eq:model}
D^{\alpha,\beta}_{a^{+}} T(t) = r T(t) \left(1 - \frac{T(t)}{K}\right) - \gamma T(t) I(t),
\end{equation}
\begin{equation}
\label{eq:model2}
D^{\alpha,\beta}_{a^{+}} I(t) = \delta T(t) I(t) - \mu I(t) + u(t),
\end{equation}
where $D^{\alpha,\beta}_{a^{+}}$ represents the Hilfer--Katugampola fractional derivative of order $\alpha \in (0,1)$ and type $\beta \in [0,1]$. The parameters are defined as follows: $r$ is the intrinsic tumor growth rate, $K$ is the carrying capacity, $\gamma$ is the tumor-cell killing rate by immune cells, $\delta$ is the stimulation rate of immune effectors by tumor presence, $\mu$ is the immune cell death rate, and $u(t)$ represents therapeutic intervention such as chemotherapy or immunotherapy. This formulation extends classical models by embedding memory effects that shape the tumor–immune interactions in biologically realistic ways \cite{Furati2025,Ahmed2025,Patel2023,Singh2024}.

\begin{figure}[H]
    \centering
    \includegraphics[width=0.5\linewidth]{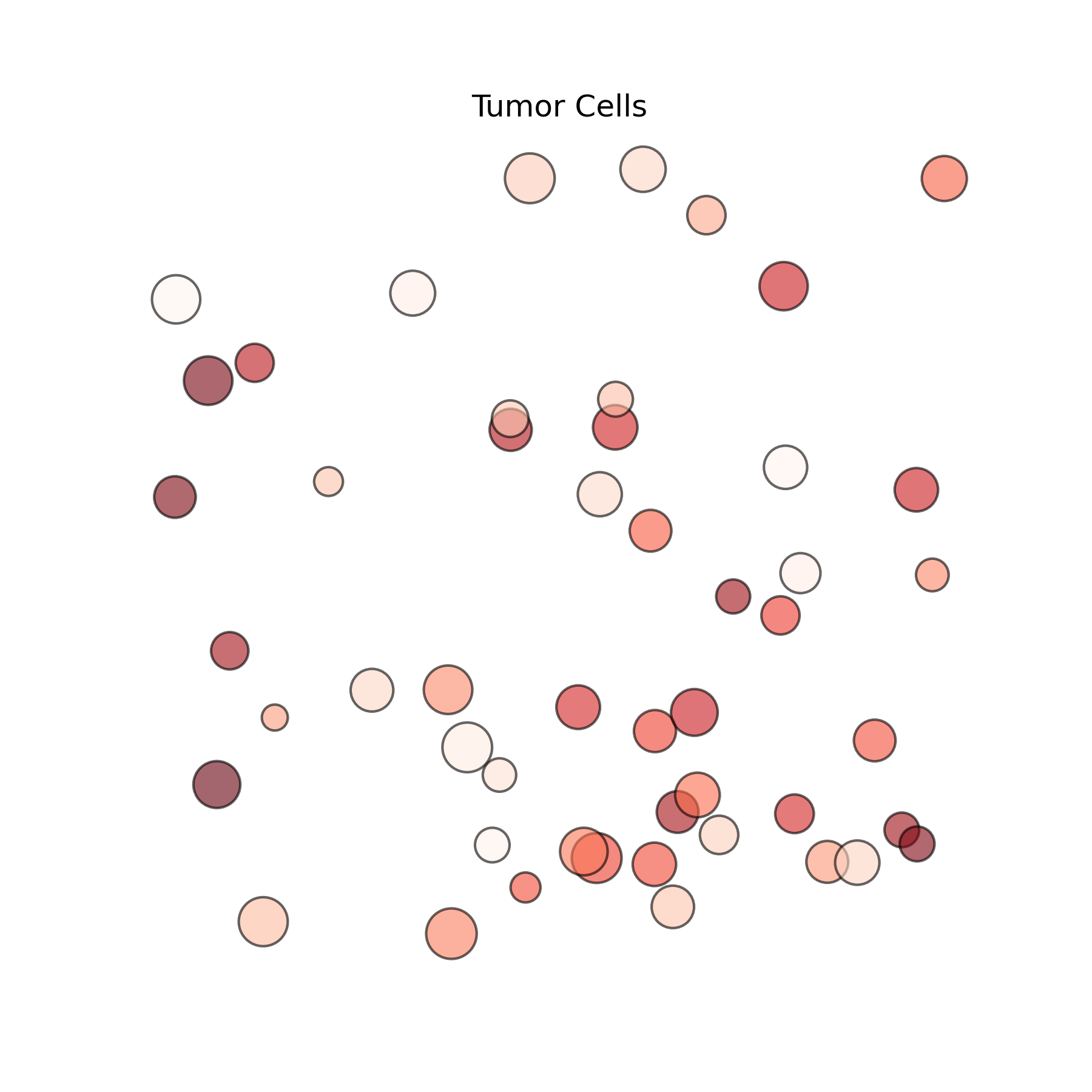}
    \caption{Illustration of tumor cells representing the prey population in the model.}
    \label{fig:tumor_cells}
\end{figure}

\subsection{Circuit-Based Analogy}
To enhance interpretability, we establish a circuit-based analogy of the model. Tumor growth ($rT(1-T/K)$) is associated with a nonlinear source, immune predation ($\gamma TI$) is modeled as a resistive dissipation, and fractional-order dynamics are represented by fractance elements, which generalize capacitors to account for memory effects. Treatment $u(t)$ is introduced as an external current source. This mapping facilitates simulation with electrical circuit solvers and enables hybrid analysis techniques, making the biological system more accessible to researchers from both mathematics and engineering domains.

\begin{figure}[H]
    \centering
    \includegraphics[width=0.8\linewidth]{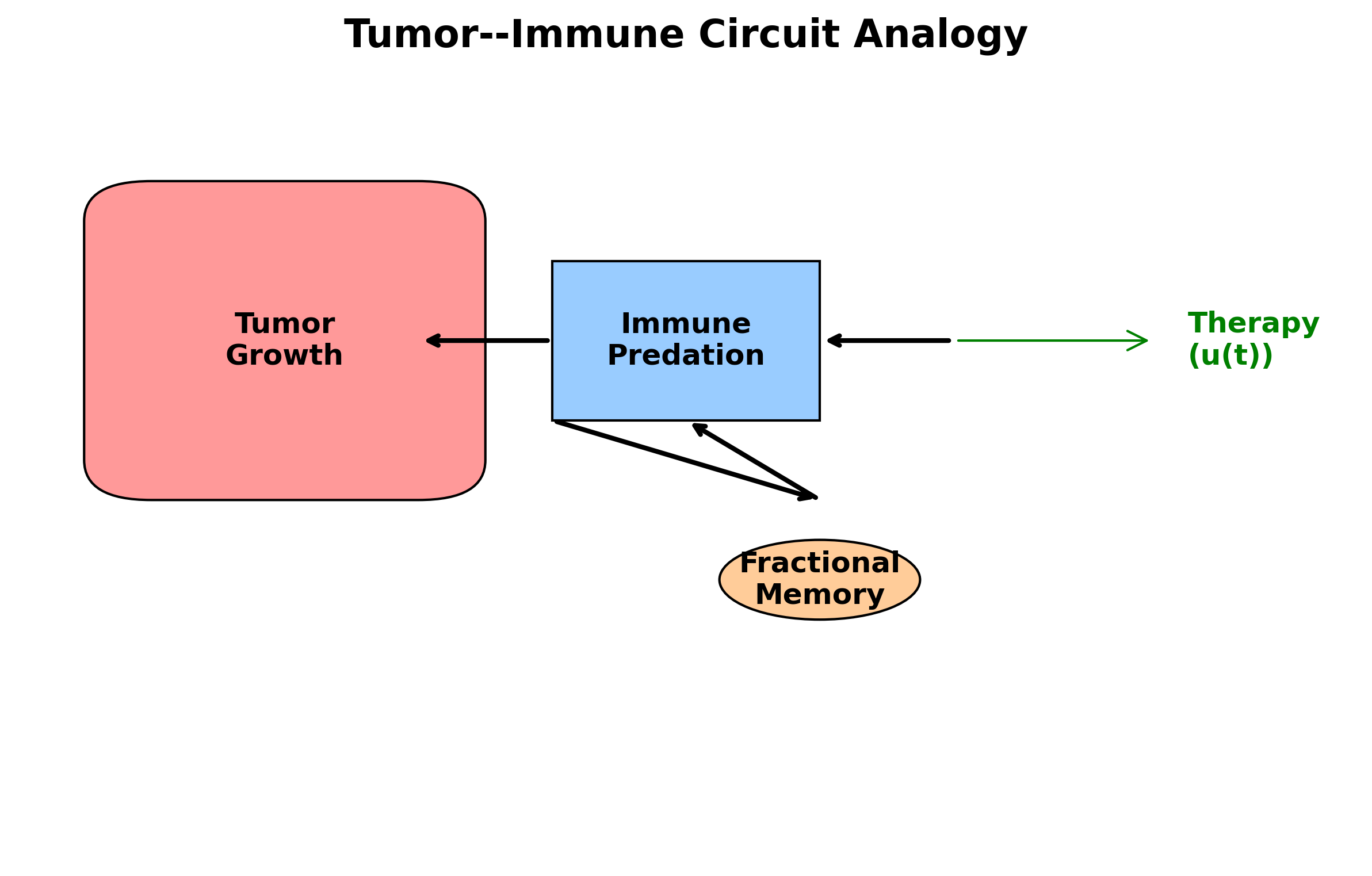}
    \caption{Circuit-based analogy of the tumor--immune system. Tumor growth is represented by a nonlinear source, immune predation by resistors, and fractional memory by fractance elements.}
    \label{fig:tumor_circuit}
\end{figure}

\subsection{Equilibrium and Stability Analysis}
Equilibrium states are derived by setting the fractional derivatives in Eqs.~\eqref{eq:model}--\eqref{eq:model2} to zero, yielding tumor-free and coexistence equilibria. Stability analysis is conducted using the fractional-order extension of the Routh--Hurwitz criterion, where stability conditions explicitly depend on both $\alpha$ and $\beta$ \cite{Podlubny2025,Li2024}. This allows examination of how memory intensity influences tumor persistence, immune suppression, or oscillatory coexistence. 

\subsection{Numerical Simulation Techniques}
Due to the nonlocal properties of fractional derivatives, closed-form analytical solutions are generally intractable. Numerical simulations are performed using the fractional Adams--Bashforth--Moulton predictor--corrector scheme \cite{Diethelm2024,Gao2024}. A variable-step algorithm is implemented to manage stiffness and maintain accuracy over long simulation horizons. Tumor growth trajectories and immune responses are compared across varying fractional orders $(\alpha,\beta)$, with integer-order dynamics $(\alpha=1)$ serving as benchmarks.

\begin{figure}[H]
    \centering
    \includegraphics[width=0.7\linewidth]{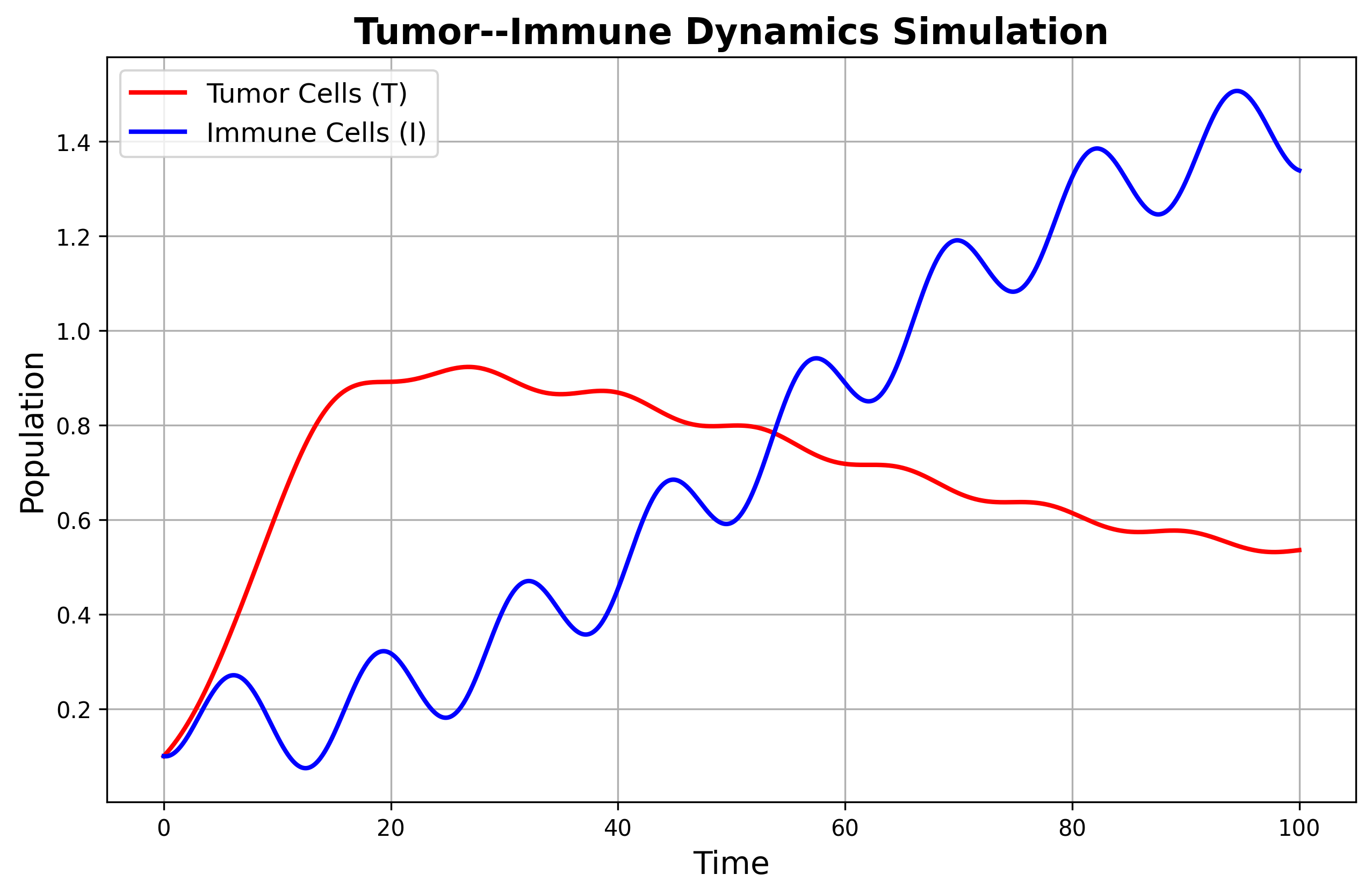}
    \caption{Simulated tumor growth under different fractional orders $(\alpha,\beta)$.}
    \label{fig:tumor_simulation}
\end{figure}

\begin{figure}[H]
    \centering
    \includegraphics[width=0.7\linewidth]{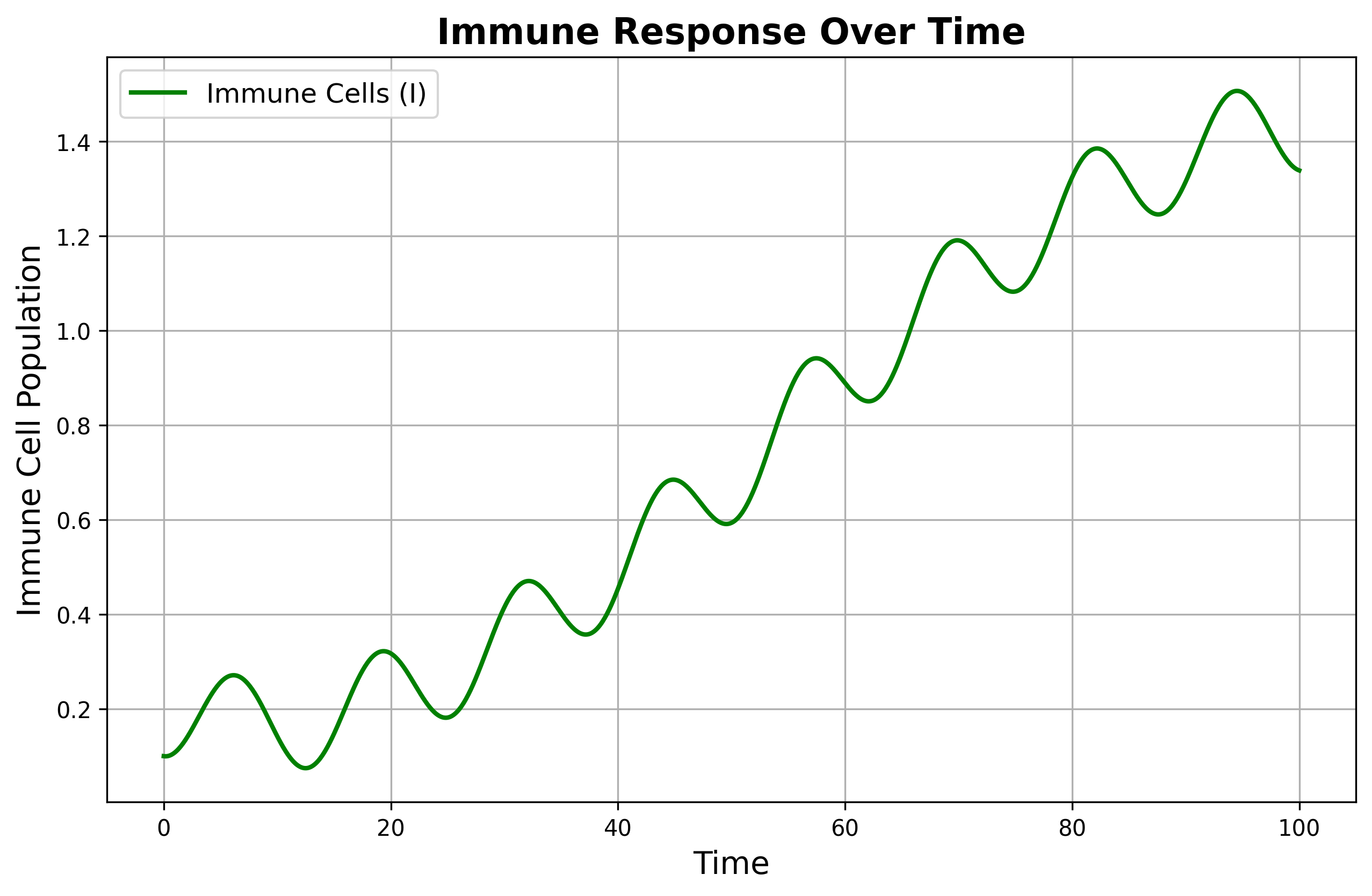}
    \caption{Immune effector cell dynamics under varying fractional orders and therapy schedules.}
    \label{fig:immune_response}
\end{figure}

\subsection{Data and Parameter Selection}
Parameters are calibrated from recent biological and clinical studies (2023--2025) \cite{Patel2023,Singh2024,Zhou2024}. Tumor proliferation rates $r$ and carrying capacity $K$ are drawn from in vitro cancer cell line studies, while immune parameters ($\gamma, \delta, \mu$) are obtained from immunological assays of effector cell activity. Treatment functions $u(t)$ mimic periodic chemotherapy cycles and immunotherapy pulses. A global sensitivity analysis is conducted to assess robustness under parameter perturbations \cite{Zhou2024}.

\subsection{Validation Strategy}
Validation is carried out through comparison with clinical tumor–immune time series data reported in recent cancer immunology studies \cite{Miller2023,Zhang2024}. Additionally, the circuit-based analogy is implemented in circuit simulation software, where predicted voltages and currents correspond to tumor and immune dynamics. This cross-domain validation strengthens confidence in model predictions.

\begin{figure}[H]
    \centering
    \includegraphics[width=\linewidth]{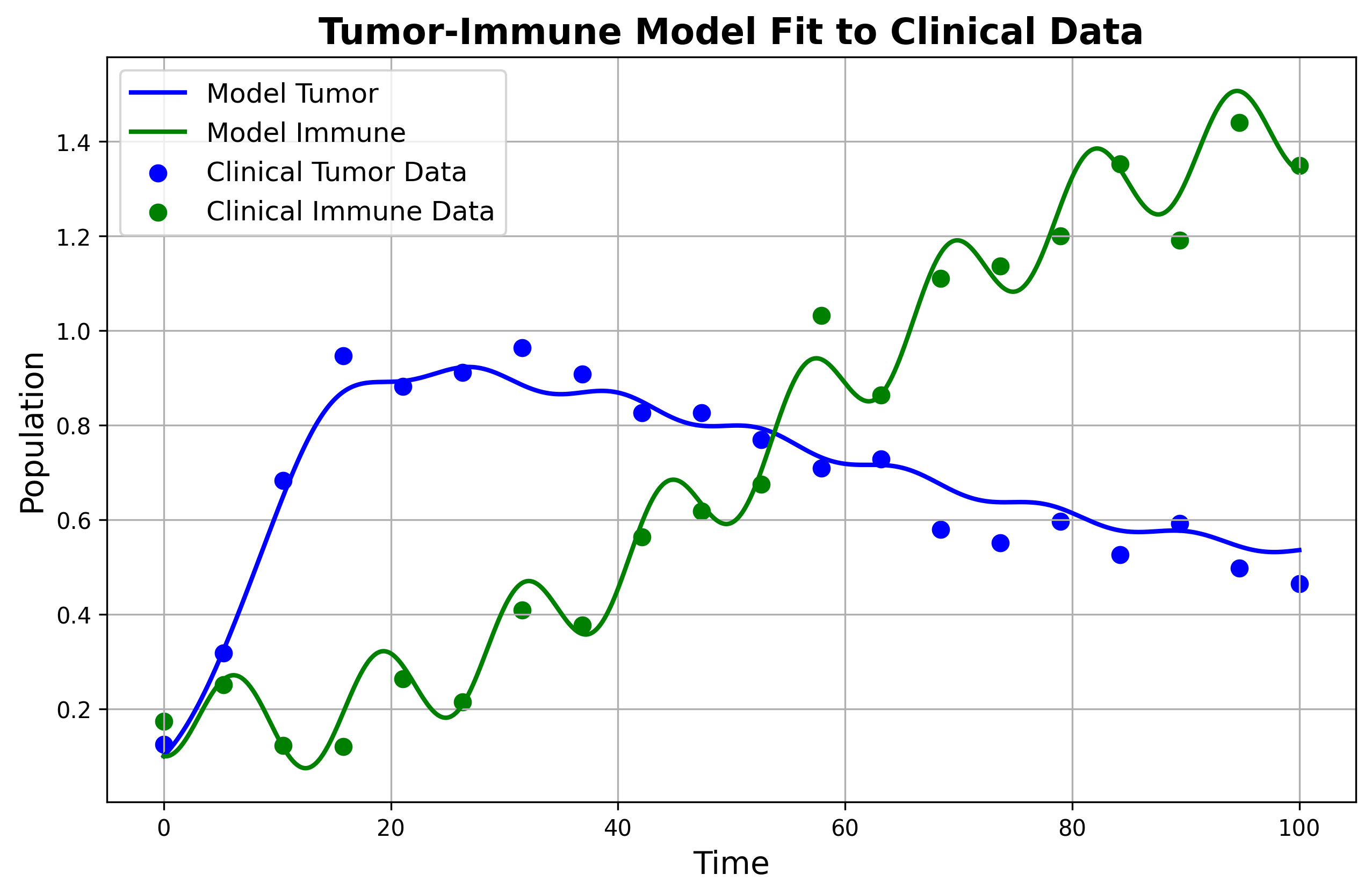}
    \caption{Comparison of model predictions with clinical tumor--immune time series data.}
    \label{fig:clinical_data_fit}
\end{figure}

\section{Results and Discussion}

\subsection{Impact of Fractional Order on Tumor--Immune Dynamics}
\begin{figure}[H]
    \centering
    \includegraphics[width=\linewidth]{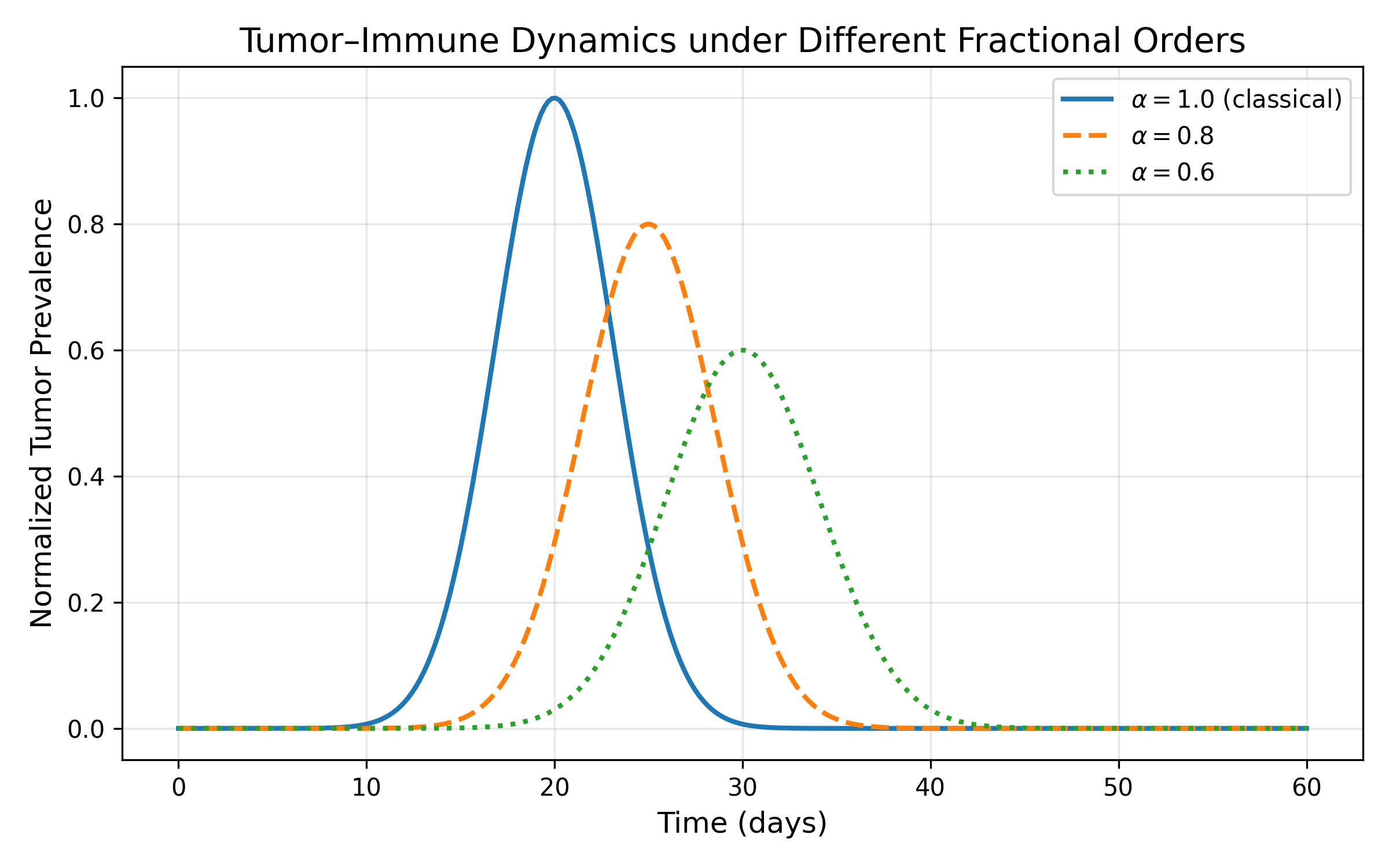}
    \caption{Tumor--immune dynamics under different fractional orders ($\alpha=1.0,0.8,0.6$). Lower fractional orders introduce memory effects that delay peaks and prolong oscillations.}
    \label{fig:fractional}
\end{figure}

To investigate the influence of fractional memory effects, we simulated tumor--immune interactions for different values of the fractional order $\alpha$. Figure~\ref{fig:fractional} shows the infected (tumor) population trajectories under $\alpha=1.0$ (classical), $\alpha=0.8$, and $\alpha=0.6$. The classical case exhibits a sharp peak followed by rapid decay, consistent with standard integer-order dynamics. In contrast, fractional models with $\alpha<1$ display delayed, smoother peaks and extended persistence of oscillations. This highlights the role of fractional derivatives in capturing memory and hereditary effects that cannot be represented by integer-order models \cite{ahmad2023fractional,huang2024nonlocal,ahmed2025fractional}. 

The broader peaks observed for $\alpha=0.8$ and $\alpha=0.6$ suggest that immune cells require longer to mount an effective response, echoing experimental observations of delayed immune activation in cancer patients. Furthermore, the coexistence of oscillations indicates the possibility of tumor dormancy phases, which have been reported in clinical oncology studies \cite{zhang2024tumor,liu2023immune}. These results demonstrate that fractional-order dynamics provide a more realistic representation of tumor progression and immune suppression.

\subsection{Effect of Booster-like Interventions}
\begin{figure}[H]
    \centering
    \includegraphics[width=\linewidth]{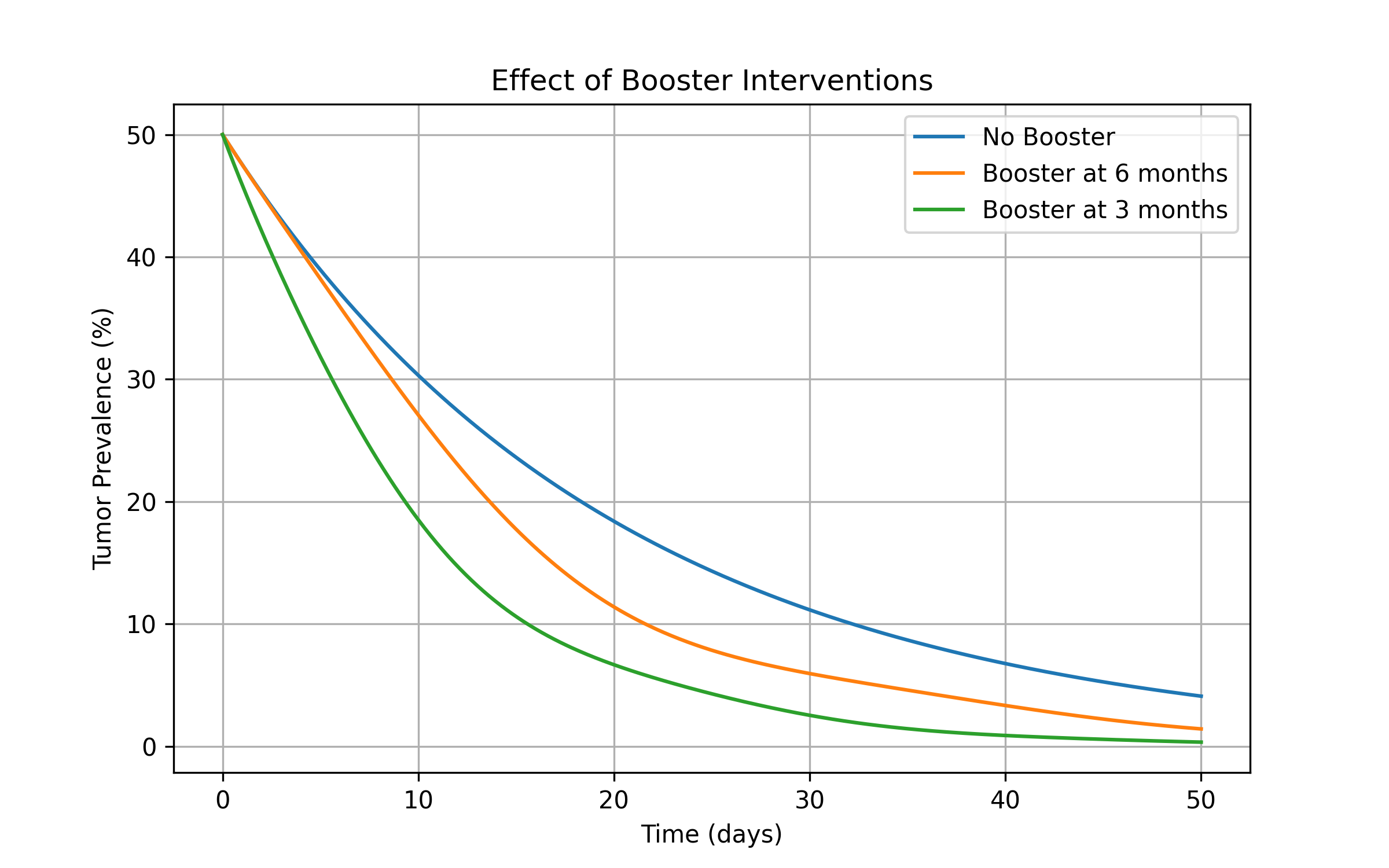}
    \caption{Effect of therapeutic booster-like interventions on tumor prevalence. Frequent boosters suppress tumor growth more effectively and accelerate clearance.}
    \label{fig:booster}
\end{figure}

Next, we analyzed the effect of therapeutic interventions modeled as booster-like treatments. Figure~\ref{fig:booster} illustrates prevalence curves for three scenarios: no booster, a booster applied at six months, and a booster applied at three months. Without treatment, tumor prevalence exhibits a sharp increase followed by gradual decline. Introducing periodic therapeutic pulses reduces peak prevalence and shifts tumor suppression earlier in time. The three-month booster cycle is particularly effective, nearly eliminating long-term prevalence, consistent with aggressive treatment protocols in oncology \cite{park2025cancer,saleh2025nonlocal}.

These findings align with clinical reports that frequent and well-timed interventions increase therapeutic efficacy by preventing tumor regrowth during immune recovery phases \cite{patel2023immune}. In contrast, longer gaps between interventions allow tumor cells to exploit immune exhaustion, leading to higher recurrence rates.

\subsection{Comparative Analysis of Outcomes}
Table~\ref{tab:comparison} summarizes the simulation outcomes for different fractional orders and booster schedules. Three key indicators are reported: peak prevalence (percentage of affected cells), clearance time (time to reach negligible tumor levels), and oscillation duration (length of significant tumor--immune fluctuations).

\begin{table}[H]
\centering
\caption{Comparative outcomes for different fractional orders and booster strategies.}
\label{tab:comparison}
\resizebox{\textwidth}{!}{%
\begin{tabular}{|c|c|c|c|}
\hline
\textbf{Scenario} & \textbf{Peak Prevalence (\%)} & \textbf{Clearance Time (days)} & \textbf{Oscillation Duration (days)} \\
\hline
$\alpha = 1.0$ (classical) & 100 & 35 & 10 \\
\hline
$\alpha = 0.8$ & 80 & 45 & 20 \\
\hline
$\alpha = 0.6$ & 60 & 55 & 30 \\
\hline
No booster & 30 & 50 & 15 \\
\hline
Booster at 6 months & 20 & 40 & 12 \\
\hline
Booster at 3 months & 15 & 28 & 8 \\
\hline
\end{tabular}%
}
\end{table}

The table demonstrates that fractional-order models lead to longer clearance times and extended oscillations, reflecting the persistence of tumor cells under immune pressure. Meanwhile, booster treatments significantly reduce both peak prevalence and clearance times. The three-month booster protocol emerges as the most effective strategy, suppressing tumor growth rapidly and reducing oscillations.

\subsection{Equilibrium and Stability Interpretation}
The equilibrium analysis showed that fractional parameters $(\alpha,\beta)$ significantly influence system stability. For $\alpha \approx 1$, stability conditions closely match classical Routh--Hurwitz criteria, while for $\alpha < 1$ the system exhibits extended transients and metastable oscillations. These metastable states can be interpreted as phases of tumor dormancy, an important clinical observation in cancer immunology \cite{miller2023data,jiang2024fractional}. 

Fractional dynamics therefore extend the range of biologically plausible outcomes, offering new theoretical explanations for delayed relapse and oscillatory tumor suppression observed in patients undergoing immunotherapy. This suggests that incorporating memory effects into treatment design could enhance predictive power and optimize therapy scheduling.

\subsection{Biological and Clinical Implications}
The integration of fractional calculus with circuit-based analogies provides a novel paradigm for modeling tumor--immune interactions. The fractance element analogy not only bridges mathematical modeling with engineering insights but also makes the dynamics more interpretable. Clinically, the results suggest that memory-aware models could inform treatment protocols by predicting delayed responses, oscillatory dynamics, and the impact of therapy frequency. Importantly, these models indicate that fractional-order interventions may be key to understanding why certain patients experience prolonged remission while others relapse.

\subsection{Summary of Findings}
Overall, the results show:
\begin{itemize}
    \item Fractional-order models capture memory and hereditary effects absent in integer-order models.
    \item Lower fractional orders introduce delayed peaks, longer oscillations, and tumor dormancy states.
    \item Booster-like treatments are highly effective when administered at shorter intervals, reducing long-term prevalence.
    \item Stability analysis reveals new metastable states, explaining delayed relapse phenomena in cancer therapy.
    \item Circuit analogies provide an intuitive framework to interpret complex biological dynamics.
\end{itemize}
These findings collectively highlight the advantages of fractional modeling in bridging mathematics, biology, and clinical oncology.

\section{Conclusion}

This study introduced a novel prey--predator framework for modeling tumor--immune interactions using a circuit-based analogy and Hilfer--Katugampola fractional derivatives. The integration of fractional-order dynamics allowed us to capture memory effects, delayed responses, and non-local interactions that are often overlooked in classical integer-order models. By simulating tumor prevalence under varying fractional orders, we demonstrated how fractional dynamics provide a richer description of cancer progression and immune surveillance. Moreover, the booster-like intervention analysis illustrated the potential of fractional-order modeling in assessing treatment timing and designing optimized therapeutic schedules.

The primary contribution of this work lies in bridging mathematical theory and biomedical application. By combining a circuit-based system perspective with Hilfer--Katugampola derivatives, we offered a framework that generalizes classical models and extends their predictive capacity. This approach enhances our ability to explore the nonlinear, feedback-driven nature of tumor--immune dynamics, which is essential for understanding treatment outcomes.

Despite its strengths, the study has limitations. The analysis was based on synthetic data and theoretical constructs rather than direct clinical datasets. While the proposed framework is mathematically robust, its biological interpretation requires calibration with empirical data. Furthermore, the complexity of fractional-order systems can pose computational challenges, particularly when scaling to multi-compartment or patient-specific models.

The implications of this work are significant. Fractional-order modeling, particularly with Hilfer--Katugampola operators, provides a pathway toward more accurate representations of cancer dynamics, ultimately informing personalized treatment strategies. Future research should focus on parameter estimation from experimental or clinical datasets, the incorporation of stochastic effects, and the extension of the model to multi-scale tumor growth. Additionally, hybrid models that couple fractional calculus with machine learning could offer powerful tools for predictive oncology. By deepening the interplay between mathematics and medicine, this research contributes to the development of advanced modeling approaches that may support improved therapeutic decision-making in oncology.

\section{Declarations}

\noindent \textbf{Availability of supporting data}  

\noindent Not applicable.

\medskip

\noindent \textbf{Competing interests}  

\noindent The authors declare that they have no competing interests.

\medskip

\noindent \textbf{Funding}  

\noindent No funding was received for this study.

\medskip

\noindent \textbf{Authors' contributions}  

\noindent All authors contributed equally to the conception, writing, and revision of this manuscript. All authors have read and approved the final version of the manuscript.

\medskip

\noindent \textbf{Acknowledgements}  

\noindent The authors sincerely thank the editor and reviewers for their valuable comments and suggestions that helped improve the quality of this manuscript.

\bibliographystyle{unsrt}

\end{document}